\documentclass[a4paper,12pt]{article}

\usepackage{graphics}
\usepackage{amsmath}
\usepackage{amsfonts}
\usepackage{amssymb}
\usepackage{latexsym}
\usepackage{color}
\usepackage{bbm}
\newcommand{\R}{\mathbb{R}}
\newcommand{\Rn}{\mathbb{R}^N}

\newcommand{\1}{\mathbbm{1}}

\usepackage{amssymb,latexsym,amsmath,amsthm,amscd,graphpap,dsfont}
\usepackage{epsfig}
\newtheorem{theorem}{Theorem}[section]
\newtheorem{mainresult}[theorem]{Main result}
\newtheorem{consequence}[theorem]{Consequences}
\newtheorem{prop}[theorem]{Proposition}

\newtheorem{lemma}[theorem]{Lemma}
\newtheorem{graph-crit}[theorem]{Graphic criterium}

\newenvironment{demo}{\noindent\textbf{Proof:}}{\hfill
$\square$\\ \medskip}

\newcommand{\Cm}{\mathcal{C}}

\newcommand{\Nm}{\mathbb{N}}
\newcommand{\Rm}{\mathbb{R}}

\newcommand{\Zm}{\mathbb{Z}}
\newcommand{\Em}{\mathcal{E}}
\newcommand{\CC}{\mathcal{C}}
\newcommand{\DD}{\mathcal{D}}

\renewcommand{\Re}{\mathrm{Re}}

\newcommand{\no}{n$^{\text{o}}$}
\newcommand{\Drond}[2]{\frac {\partial #1}{\partial #2}}

\numberwithin{equation}{section}

\newcounter{constantes}

\newdimen\texpscorrection
\newdimen\figcenter
\def\figurewithtex #1 #2 #3 #4 #5\cr{\null
  {\goodbreak\figcenter=\hsize\relax
  \advance\figcenter by -#4truecm
  \divide\figcenter by 2
  \begin{figure}[hbt]
  \vskip #3truecm\noindent\hskip\figcenter
  \includegraphics{#1}{\hskip\texpscorrection\input #2 }
  \vskip 0.8truecm{\baselineskip=0.8\baselineskip
  \noindent \vbox{\noindent {\footnotesize #5}}\par}
  \end{figure}}}
\def\point#1 #2 #3 {\rlap{\kern #1 truecm
\raise #2 truecm \hbox{#3}}}

\begin{document}

\title{\bf Asymptotic profiles for a travelling front solution of a biological equation } 

\author{ Guillemette CHAPUISAT\footnote{Universit\'e Paul C\'ezanne,
    CNRS, LATP (UMR 6632), Facult\'e des Sciences et Techniques de St
    J\'er\^ome, Case Cour A, av Escadrille Normandie-Niemen, 13397
    Marseille Cedex 20, France. {\rm
      guillemette.chapuisat@univ-cezanne.fr}}~~\&~Romain 
JOLY\footnote{Institut Fourier, UMR 5582 CNRS/Universit\'e de Grenoble, 100, rue
des maths, BP 74, 38402 Saint-Martin-d'H\`eres, France. {\rm
romain.joly@ujf-grenoble.fr}
}}

\maketitle

\begin{abstract}
We are interested in the existence of depolarization waves in the human brain.
These waves propagate in the grey matter and are absorbed in the white matter.
We
consider a two-dimensional model  $u_t=\Delta u + f(u) \1_{|y|\leq R} -
\alpha u \1_{|y|>R}$, with $f$ a bistable nonlinearity taking effect only
on the domain $\Rm\times [-R,R]$, which represents the grey matter
layer. We study the existence,
the stability and the energy of non-trivial asymptotic profiles of
possible travelling fronts. For this purpose, we present dynamical systems
technics and graphic criteria based on Sturm-Liouville theory and apply them to
the above equation. This yields three different behaviours of the
solution $u$ after stimulation, depending of the thickness $R$ of the grey
matter. This may partly explain the difficulties to observe depolarization
waves in the human brain and the failure of several therapeutic trials.\\[2mm]
{\bf AMS classification codes (2000):} 34C10, 35B35, 35K57, 92C20.
\end{abstract}

\section{Introduction}\label{intro}
The propagation of depolarization waves, also called spreading depressions,
may appear in a brain during
strokes, migraines with aura or epilepsy. In rodent brains,
the propagation of depolarization waves has been observed for more than fifty
years \cite{Leao}. During stroke in rodent brain, they cause important
damages and are therefore a therapeutic target. Pharmacological agents
blocking the appearance of 
those waves have been studied and reduce strongly after-effects of
stroke in the rodent brain \cite{Mies, Obeidat}. However propagation
of depolarization waves during stroke in the human brain is still
uncertain \cite{aitken, mayevsky, maclachlan,somjen,sramka} and the
pharmalogical agents used in rodent have seemed to have no effect on
human stroke. Mathematical models of these depolarization waves may
help to understand these points.

In \cite{GTN}, the first author has build a mathematical model of these waves
and has studied numerically how the morphology of the brain may influence
their propagation. Simplifications of the mathematical model lead to the study
of propagation phenomena for the following biological equation 
\begin{equation}\label{eq}
\frac{\partial u}{\partial t} - \triangle u = f(u) \1_{\Omega} - 
                                    \alpha u \1_{\Rn \setminus \Omega},  
\qquad t \in \R, \, X \in \Rn,
\end{equation}
\noindent where $f(u)=\lambda \; u (u- a) (1-u)$ is the usual bistable
nonlinearity with $a \in ]0,1/2[$ and $\lambda>0$, and where
$\alpha$ is a positive number. We denote by $\1_\Omega$ the 
caracteristic function of the domain $\Omega$ that is
$\1_\Omega(X)=1$ if $X\in\Omega$ and $\1_\Omega(X)=0$ elsewhere. The
function $u$ represents the depolarization of the brain that is that
if $u(X)=0$ the brain is normally polarized at the point $X$ whereas if $u(X)=1$
the brain is totally depolarized. The domain $\Omega$
represents the grey matter of a brain, where the
reaction-diffusion process that triggers depolarization waves takes
place, whereas $\Rm ^N \setminus\Omega$ contains the white matter,
where the waves are absorbed. The problem is to understand the
influence of the geometry of $\Omega$ on the propagation of waves. In
particular, the layer of grey matter in the human brain is thiner than
in other species and admits more circumvolutions. The influence of 
the circumvolutions of the grey matter is partially studied in \cite{CPDE}. 

We will focus here on the
influence of the thickness of the layer of grey matter. We may assume that
$\Omega$ is a straight cylinder of radius $2R$. In \cite{GTN}, numerical studies
have shown that small values of $R$ may prevent the spreading of the ionic
waves. A partial theoretical study in any dimension $N$ has been led in
\cite{GC}, where the first author has proved the non-existence of the
depolarization waves if the thickness $R$ is small enough and their existence if
$R$ is large enough.

The results of \cite{GC} are not completely satisfactory. First,
they only deal with $R$ small or large enough. Secondly, no complete
description of the possible asymptotic profiles of the travelling waves, their
stability, or their energy, has
been pursued. As a consequence, even if the existence of waves is proved,
nothing is known about their asymptotic profiles and their stability. The
purpose of the present paper is to complete this study for any $R$, in
dimension $N=2$. More precisely, we consider the equation
\begin{equation}\label{eq-2d}
\frac{\partial u}{\partial t}=\triangle u + f(u) \1_{|y|\leq R} - 
                                    \alpha u \1_{|y|>R}, \qquad t \in \R, \,
(x,y)\in \Rm^2~.
\end{equation}
We are looking to solutions travelling in the $x$ direction at speed
$c$, that are solutions $u$ of \eqref{eq-2d} which can be written
$u(x,y,t)=v(x-ct,y)$. Travelling fronts are solutions
$u(x,y,t)=v(x-ct,y)$ such that there are two asymptotic profiles
$V_-\neq V_+$ with $\lim_{\xi\rightarrow \pm \infty}
\|v(\xi,.)-V_\pm\|_{H^1(\Rm)} =0$. Using standard elliptic
estimates, a profile $V$ is a solution in $H^1(\Rm)$ of the equation 
\begin{equation}\label{eq-profil} 
V''+f(V)\1_{|y|\leq R}-\alpha V\1_{|y|>R}=0,~~~y\in\Rm~.
\end{equation} 
Notice that a profile is trivially associated to an equilibrium point
$E(x,y)=V(y)$ that is a stationnary solution of \eqref{eq-2d}. The trivial
equilibrium point $E\equiv 0$ corresponds to the normal state of the brain,
whereas a non-trivial equilibrium point corresponds to a depolarized state,
deleterious for the brain. We are more particularly interested in the existence
of travelling fronts with positive speed $c$ connecting the profile $V_+\equiv
0$ with a non-trivial profile $V_-$. Such fronts correspond to the invasion of
the equilibrium state $0$ by a deleterious state.

The purpose of this article is to achieve a complete study of the
existence of the profiles, their stability and their energy (see the
following sections for the definitions of stability and energy).
\begin{mainresult}\label{th-1}
There exist two critical thicknesses $R_1>R_0>0$ such that:\\
i) if $0<R<R_0$, there is no non-trivial profile, i.e. that $V\equiv 0$ is
the only solution of \eqref{eq-profil} in $H^1(\Rm)$.\\
ii) If $R_0<R$, there exist non-trivial profiles. One of them, denoted
by $V_M$, is larger than every other one. The largest profile $V_M$ is
stable, and every other non-trivial profiles are
unstable.\\
iii) The energy of the unstable profiles is always larger than the energy of
the stable profiles $0$ and $V_M$. If $R_0<R<R_1$, the energy of $V_M$ is
larger than the energy of $0$, whereas it is smaller if $R>R_1$.
\end{mainresult}

In this paper, we first adapt several
graphic criteria for stability of profiles, based on Sturm-Liouville
theory. Then, we present dynamical systems technics for studying energy of
profiles. These
criteria and technics are general ones and are valid for any nonlinearity
$f$. We obtain the above results by applying these ideas to the particular
nonlinearity $f(u)=\lambda \; u (u- a) (1-u)$. The check on the graphic criteria
sometimes relies on a numerical study of the phase plane.

Once the profiles and their stability completely described by the Main
Result \ref{th-1}, one is able to obtain the existence of travelling fronts.
\begin{consequence}\label{th-2}
$~$\\
i) If $0<R<R_0$ then there is no travelling front for Equation \eqref{eq-2d}.\\
ii) If $R_0<R<R_1$ there exists a globally stable travelling front
$u(x,y,t)=v(x-ct,y)$ with $\lim_{\xi\rightarrow +\infty} v(\xi,.)=0$ and
$\lim_{\xi\rightarrow -\infty} v(\xi,.)=V_M$ with a negative speed $c<0$.\\
ii) If $R_1<R$ there exists a globally stable travelling front
$u(x,y,t)=v(x-ct,y)$ with $\lim_{\xi\rightarrow +\infty} v(\xi,.)=0$ and
$\lim_{\xi\rightarrow -\infty} v(\xi,.)=V_M$ with a positive speed $c>0$.\\
\end{consequence}
Assertion i) is a clear consequence of the Main Result \ref{th-1}. The proof
that this Main Result also implies the other assertions is not the
goal of this paper 
and will not be detailled. One way to prove it is to use a method recently
introduced by Risler in \cite{Risler}. The fundamental idea is to use the
existence of an energy functional for \eqref{eq-2d} in any Galilean frame
travelling at constant speed $s$ in the $x$ direction. Keeping this basic idea,
the original proof has been revised in \cite{Gallay-Risler} and
\cite{Gallay-Joly}. The adaptation of these technics to Equation \eqref{eq-2d}
has been detailled in \cite{GC}. The original motivation of Risler's technics
was to develop a method of proof of existence and stability of fronts in
equations where no comparison principle is available. As suggested by the work
of the first author in \cite{GC}, the method of Risler is also usefull to get
the existence and stability of fronts in equations as \eqref{eq-2d}
for which even if there is a comparison principle, no particular positive
subsolution is known.

\begin{figure}[htp]
\begin{center}
\epsfig{file=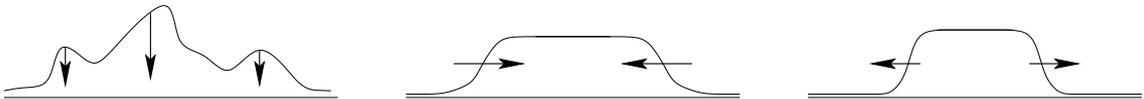,width=15cm}
\caption{\it The three typical behaviours of the depolarization in the grey
matter
after stimulation. From left to right: $R<R_0$, $R_0<R<R_1$ and $R_1<R$.}
\label{fig-intro}
\end{center}
\end{figure}

The results of this article yield a good insight into the different behaviours
of the depolarization in the grey matter after it has been stimulated. If the
stimulation takes place in a part of a brain where the grey matter is thin
($R<R_0$), the neurons of the grey matter quickly repolarize, i.e. the
solution of \eqref{eq-2d} goes back to zero, uniformly and exponentially fast.
If the grey matter is slightly thicker ($R_0<R<R_1$), the repolarization is
slower and more progressive, the depolarized area disappearing by shrinking.
Mathematically speaking, there exists a stable non-trivial profile but with an
energy larger than zero. Thus, the equilibrium state zero invades the excited
state by travelling fronts, reducing the excited area. Notice that it may take
very long time to go back to rest if the initial excited area is large.
Finally, if the grey matter is thick ($R_1<R$), a depolarization wave
propagates. Indeed, the excited state has a lower energy than zero.
These three behaviours are illustrated in Figure \ref{fig-intro}. The fact that
the behaviour depends on the width of the grey matter may explain why the
 depolarization waves have been observed or not in the human brain depending on
the experiments, see \cite{aitken, mayevsky, maclachlan} or\cite{sramka} and
the discussion of Section \ref{sect-discussion}. 

The paper is splitted as follows. In Section \ref{sect-para} we explicit a
relation between the profiles and the equilibrium points of a parabolic equation
in $(-R,R)$. In Section \ref{sect-exist}, we obtain graphic criteria of
existence and stability of the profiles, mainly based on Sturm-Liouville theory.
The energy of the different profiles is  studied in Section \ref{sect-energy},
using technics coming from the infinite dimensional dynamical systems theory.
Finally, in Section \ref{sect-discussion}, we discuss the relations between our
mathematical results and the biological phenomena of spreading waves in
the brains.

%%%%%%%%%%%%%%%%%%%%%%%%%%%%%%%%%%%%%%%%%%%%%%%%%%%%%%%%%%%%%%%%
%%%%%%%%%%%%%%%%%%%%%%%%%%%%%%%%%%%%%%%%%%%%%%%%%%%%%%%%%%%%%%%%
%%%%%%%%%%%%%%%%%%%%%%%%%%%%%%%%%%%%%%%%%%%%%%%%%%%%%%%%%%%%%%%%

\section{Relations with a parabolic equation on $(-R,R)$}\label{sect-para}
In this preliminary section, we enhance an obvious relation between the
profiles and the equilibrium points of a parabolic equation on the segment
$(-R,R)$. We also recall some basic properties of the dynamics of this parabolic
equation. In fact, in this paper, we could perform all the arguments with the
parabolic equation on the whole domain $\Rm$. However, this would bring useless 
technicalities.

\subsection{Back to a bounded interval}
The profiles $V$ are the solutions of Equation \eqref{eq-profil}. In other
word, they are the stationnary solutions of the evolution equation 
\begin{equation}\label{eq-evo-prof} 
\Drond ut=\Drond {^2u}{y^2}+ g(y,u)
\end{equation}
where $g(y,u)=f(u)\1_{|y|\leq R}-\alpha u \1_{|y|>R}$. The energy corresponding
to this reaction-diffusion equation is the functional $E:H^1(\Rm)\longrightarrow
\Rm$ defined by 
\begin{equation}\label{eq-energy}
E(u)=\frac12 \int_\Rm \left|\frac{du}{dy}(y)\right|^2dy~-~\int_\Rm G(y,u(y))dy 
\end{equation}
where $G(y,u)=\int_0^u g(y,v)dv$ is a primitive of $g$.\\[1mm]

We notice that any profile $V$ satisfies $V(y)=c_+ e^{\sqrt\alpha y}+ c_-
e^{-\sqrt \alpha y}$ outside $(-R,R)$. Since $V$
belongs to $H^2(\Rm)$, one has $c_+=0$ for $y>R$ and $c_-=0$ for
$y<R$. Since $H^2(\Rm)\subset \Cm^1(\Rm)$, $V$ is continuous and solving
\eqref{eq-profil} is equivalent to solving the following problem 
\begin{equation}\label{eq-profile-2}
\left\{\begin{array}{ll} V''(y)+f(V(y))=0 & ~~y\in(-R,R) \\
V'(-R)=\sqrt{\alpha}\; V(-R) & \\  V'(R)=-\sqrt{\alpha}\; V(R)~. &
\end{array}\right.
\end{equation}
Equation \eqref{eq-profile-2} caracterises the equilibrium points of a
parabolic equation on $(-R,R)$ with Robin boundary conditions
\begin{equation}\label{eq-evo-prof-2} 
\left\{\begin{array}{ll} 
\displaystyle \Drond{u}{t}(y,t)=\Drond{^2u}{y^2}(y,t)+f(u(y,t)) &
~~(y,t)\in(-R,R)\times\Rm_+\\ 
\displaystyle \Drond uy (-R,t)=\sqrt{\alpha}\; u(-R,t) &~~t>0\\  
\displaystyle \Drond uy
(R,t)=-\sqrt{\alpha}\; u(R,t) &~~t>0
\end{array}\right.
\end{equation}
The energy corresponding to Equation \eqref{eq-evo-prof-2} is given by
\begin{equation}\label{eq-energy-2}
\Em(u)=\frac12 \int_{-R}^R \left|\frac{du}{dy}(y)\right|^2dy~-~\int_{-R}^R 
F(u(y))dy + \frac {\sqrt \alpha}2 \left(|u(R)|^2+|u(-R)|^2\right)
\end{equation}
where $F(u)=\int_0^u f(v)dv$ is a primitive of $f$.

It is straightforward to verify that if $V$ is a profile, then
$E(V)=\Em(V_{|(-R,R)})$. As a conclusion, the problems of the existence and the
energy of the profiles
are equivalent to the study of the equilibrium points of the parabolic equation
\eqref{eq-evo-prof-2} and their energy. We will often abusively denote
$V_{|(R,R)}$ by simply $V$. However, notice that the dynamics of the equation
\eqref{eq-evo-prof} on the whole real line and the ones of \eqref{eq-evo-prof-2}
on $(-R,R)$ are different. Moreover, the spectrum of the linearization of
\eqref{eq-evo-prof} at a profile $V$ is different from the spectrum of the
linearization of \eqref{eq-evo-prof-2} at $V_{|(-R,R)}$, even if we prove in
Section \ref{subsection-exist-vp} that the number of positive eigenvalues is
the same.

\subsection{Basic facts about the parabolic equation on $(-R,R)$}

The theory of the parabolic equation \eqref{eq-evo-prof-2} with Robin boundary
conditions is very well-known, see for example \cite{Henry}. We recall here
some basic dynamical properties, which will be used in this paper.

First, \eqref{eq-evo-prof-2} admits a comparison principle.
\begin{lemma}\label{lemme-compare}
Let $u(t)$ and $v(t)$ be two solutions of \eqref{eq-evo-prof-2} with
respective initial data $u_0$ and $v_0$ belonging to $H^1(-R,R)$. Assume
that $u_0(y)\leq v_0(y)$ for all $y\in (-R,R)$. Then, for all
$(y,t)\in (-R,R) \times\Rm_+$, $u(y,t)\leq v(y,t)$. 
\end{lemma}
Secondly, the parabolic type of \eqref{eq-evo-prof-2} and the fact that
$f(u)u<0$ for large $u$ implies that any solution $u(t)$ with initial data
$u_0\in H^1(-R,R)$ is bounded in $H^k(-R,R)$ uniformly for $t\in
[t_0,+\infty)$ for any $t_0>0$ and for any $k\in\Nm$. This yields a compacity
result.
\begin{lemma}\label{lemme-compacite}
Let $(t_n)_{n\in\Nm}$ be a sequence of times such that
$\liminf t_n>0$, let $(u^0_n)_{n\in\Nm}\subset H^1(-R,R)$ be a sequence of
initial data and let $u_n(t)$ be the solutions of \eqref{eq-evo-prof-2} with
$u_n(0)=u_n^0$.\\
Then, there exist an extraction $\varphi$ and a function $u_\infty\in H^1(-R,R)$
such that
$\|u_{\varphi(n)}(t_{\varphi(n)})-u_\infty\|_{H^1(-R,R)}\longrightarrow 0$
when $n$ goes to $+\infty$.
\end{lemma}
The energy $\Em$ defined by
\eqref{eq-energy-2} is a strict Lyapounov functional for the dynamics of
\eqref{eq-evo-prof-2}: for
each $u(t)\in H^1(-R,R)$ solution of \eqref{eq-evo-prof-2}, the function
$\Em(u(t))$ is decreasing in time, except if $u(t)$ is constant, that is except
if $u$ is an equilibrium point, and thus a profile. Indeed, it is
straightforward to get $\partial_t \Em(u(t))=-\int_{-R}^R |\Drond ut (y,t)|^2
dt$. The existence of a strict Lyapounov function yields what is called
gradient dynamics. In addition of the compacity property of Lemma
\ref{lemme-compacite}, it is very classical to obtain Lasalle
principle, see \cite{Hale}. 
\begin{lemma}\label{lemme-gradient}
Let $u(t)$, $t\geq0$, be a solution in $H^1(-R,R)$ of \eqref{eq-evo-prof-2}. Let
$\omega(u)$ be the $\omega-$limit set of $u$, that is 
$$\omega(u)=\left\{ v\in H^1(-R,R)~/~\exists
(t_n)_{n\in\Nm},~t_n\xrightarrow[n\rightarrow +\infty]{}+\infty~\text{
  and }~\|u(t_n)-v\|_{H^1(-R,R)}\xrightarrow[n\rightarrow +\infty]{} 0
\right\}~. $$ 
Then, $\omega(u)$ is a non-empty connected compact set which only consists in
equilibrium points, i.e. in profiles.\\ 
Similarly, if $u(t)$ is defined for all $t\in\Rm$ and uniformly
bounded, the $\alpha-$limit set  
$$\alpha(u)=\left\{ v\in H^1(-R,R)~/~\exists
(t_n)_{n\in\Nm},~t_n\xrightarrow[n\rightarrow +\infty]{}-\infty~\text{
  and }~\|u(t_n)-v\|_{H^1(-R,R)}\xrightarrow[n\rightarrow +\infty]{} 0
\right\}~ $$ 
is a non-empty connected compact set which only consists in
profiles. Moreover, if $V_-\in\alpha(u)$ and $V_+\in\omega(u)$ then
$\Em(V_-)>\Em(V_+)$ and in particular $V_-\not\equiv V_+$.
\end{lemma}

%%%%%%%%%%%%%%%%%%%%%%%%%%%%%%%%%%%%%%%%%%%%%%%%%%%%%%%%%%%%%%%%
%%%%%%%%%%%%%%%%%%%%%%%%%%%%%%%%%%%%%%%%%%%%%%%%%%%%%%%%%%%%%%%%
%%%%%%%%%%%%%%%%%%%%%%%%%%%%%%%%%%%%%%%%%%%%%%%%%%%%%%%%%%%%%%%%

\section{Existence and stability of profiles}\label{sect-exist}
We are interested in the asymptotic profiles of the possible travelling fronts 
of Equation \eqref{eq-2d}. We recall that a profile $V\in H^1(\Rm)$ is
solution of \eqref{eq-profil}. We wonder if there exist non-trivial
profiles (i.e. profiles $V\not\equiv 0$) and if they are stable. A profile $V$
is stable if the linear operator
$L_V: H^2(\Rm)\longrightarrow L^2(\Rm)$ defined by 
\begin{equation}\label{eq-LV}
L_V \varphi=\varphi''+\left(f'(V)\1_{|y|\leq R}-\alpha\1_{|y|>R}\right)\varphi
\end{equation}
has no spectrum in the half-plane $\{z\in\Cm, Re(z)\geq 0\}$.\\

In this section, we give graphic criteria for
existence and stability of profiles and we apply them to Equation
\eqref{eq-profil} to obtain the following result.
\begin{prop}\label{prop-profiles}
For any given $\alpha$, $\lambda$ and $a$, there exists a
positive number $R_0$ such that:\\ 
i) if $R<R_0$, then there does not exist any profile different from $0$.\\
ii) if $R\geq R_0$, then there exist non-trivial profiles. One of
them, denoted by $V_M$ is strictly larger than every others, smaller
than $1$ and is even. Moreover, if $R>R_0$ then $V_M$ and $0$ are
stable profiles and are the only ones to be
stable.\\
iii) if $\alpha\geq  a\lambda$ then if $R>R_0$, there exist
exactly two non-trivial profiles.\\  
\end{prop}

%%%%%%%%%%%%%%%%%%%%%%%%%%%%%%%%%%%%%%%%%%%%%%%%%%%%%%%%%%%%%%%%

\subsection{A graphic criterium for existence of profiles}\label{subsection-exist-profile}
We recall here the simple graphic arguments already mentionned in
\cite{GC}. As already noticed, looking for profiles $V$
solutions of \eqref{eq-profil} is equivalent to looking for solutions $V$ of
the equation \eqref{eq-profile-2}. In the phase plane, the flow
$\Phi_{y_0}:\Rm^2\longrightarrow \Rm^2$ associated to the differential equation
$V''+f(V)=0$ between $-R$ and $y_0$ is given by
$$\Phi_{y_0}(V_0,V_1)=(V(y_0),V'(y_0)) ~~\text{ where
}V \text{ solves }~~\left\{\begin{array}{l}
V''(y)+f(V(y))=0~,~~y\in (-R,y_0)\\
(V(-R),V'(-R))=(V_0,V_1)\end{array}\right. $$
Let $\DD$ and $\DD'\subset\Rm^2$ be the straight lines
$\Rm.(1,\sqrt{\alpha})$ and $\Rm.(1,-\sqrt{\alpha})$. Notice that the
conditions $V'(-R)=\sqrt{\alpha}\;V(-R)$ and  $V'(R)=-\sqrt{\alpha}\;
V(R)$ correspond to $(V,V')(-R)\in \DD$ and $(V,V')(R)\in \DD'$
respectively. Hence, the graphic interpretation of
\eqref{eq-profile-2} is the following.
\begin{graph-crit}\label{gc1}
The profiles, i.e. the solutions of \eqref{eq-profile-2}, correspond
to the trajectories of the flow $\Phi_y$ that start on $\DD$ when $y=-R$ and
finish on $\DD'$ when $y=R$. Therefore, they correspond to the
intersections of the curve $\Phi_{R}\DD$ with the straight line $\DD'$.
\end{graph-crit}

\begin{figure}[htp]
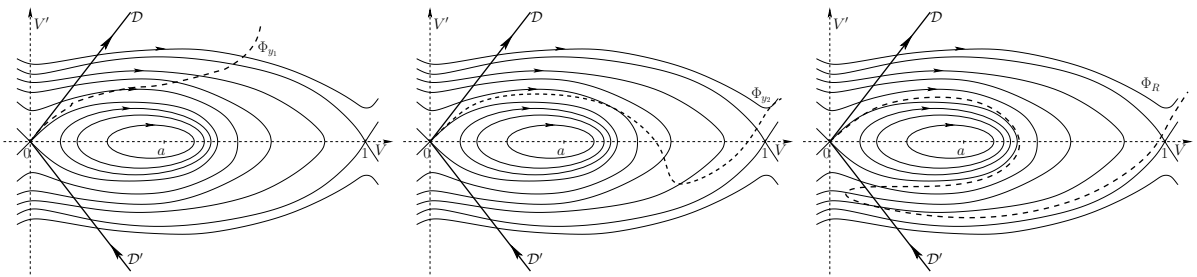

\begin{center}
\resizebox{0.32\textwidth}{!}{\input{phiR1.pstex_t}}
\resizebox{0.32\textwidth}{!}{\input{phiR2.pstex_t}}
\resizebox{0.32\textwidth}{!}{\input{phiR3.pstex_t}}
\caption{\it The phase plane for $\alpha \geq  a\lambda$ and $\Phi_R\DD$ for
three
increasing values of $R$: two cases with no non-trivial
profiles and the last with two non-trivial profiles.}
\label{fig-R1}
\end{center}
\end{figure}

\begin{figure}[htp]
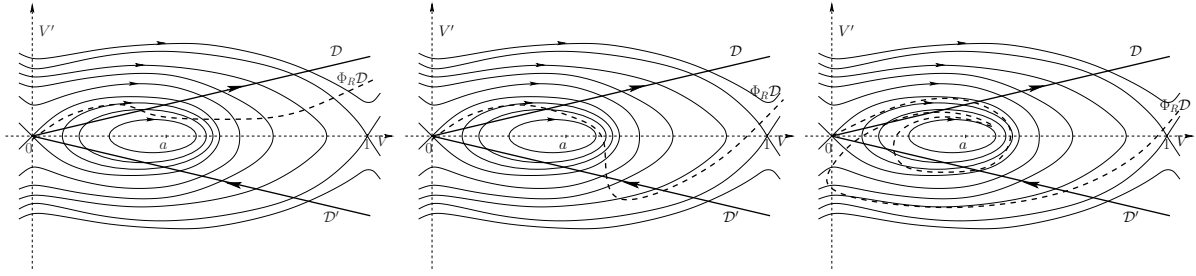

\begin{center}
\resizebox{0.32\textwidth}{!}{\input{phiR4.pstex_t}}
\resizebox{0.32\textwidth}{!}{\input{phiR5.pstex_t}}
\resizebox{0.32\textwidth}{!}{\input{phiR6.pstex_t}}
\caption{\it The phase plane for $\alpha < a\lambda$ and $\Phi_R\DD$ for three
increasing values of $R$: no non-trivial profiles, two ones and four ones.}
\label{fig-R2}
\end{center}
\end{figure}

\begin{figure}[htp]
\begin{center}
\resizebox{0.45\textwidth}{!}{\input{phiR56.pstex_t}}
\resizebox{0.4\textwidth}{!}{\epsfig{file=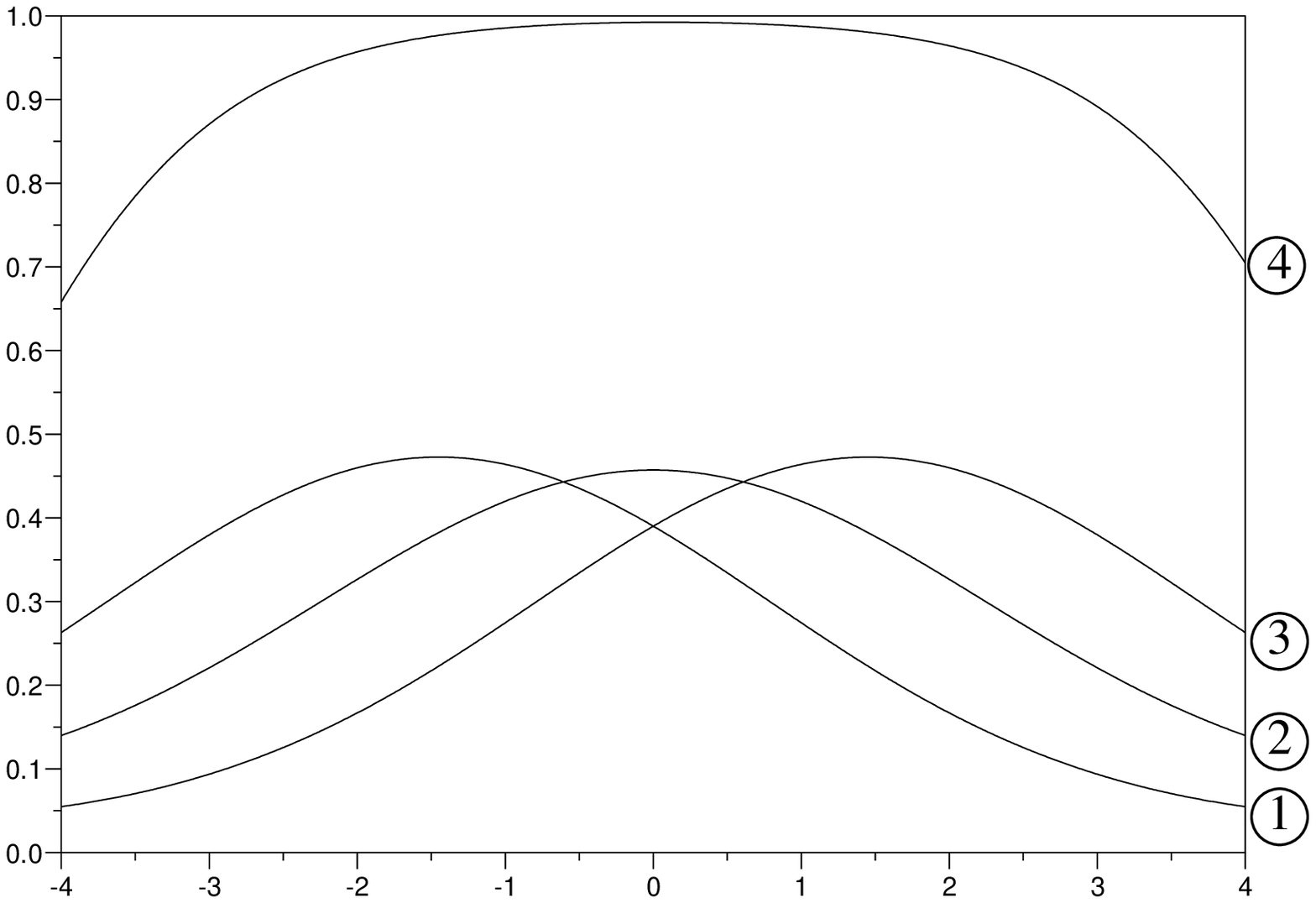, width=6cm}}
\caption{\it The intersections of $\Phi_R\DD$ and $\DD'$ in the phase plane and
the
corresponding profiles in the $(y,V)$ plane}
\label{fig-R3}
\end{center}
\end{figure}

The results of existence stated in Proposition \ref{prop-profiles}
follow from Graphic Criterium \ref{gc1}. The critical radius $R_0$ is
the first (positive) $R$ such that $\Phi_{R}\DD$ intersects $\DD'$
elsewhere than at $(0,0)$. In Figures \ref{fig-R1} and \ref{fig-R2}, we
illustrate the phase plane, the flow $\Phi_y$ and the curve $\Phi_R\DD$ for
different cases. Figure \ref{fig-R3} shows the graph of four profiles $V(y)$ in
$[-R,R]$ and the corresponding intersections in the phase plane.


Figures \ref{fig-R1} and \ref{fig-R2} show that there are two cases. If $\DD$
is always outside the domain delimited by the homoclinic orbit to zero, then
there cannot exist more than two non-trivial profiles. If $\DD$ intersects the
homoclinic orbit to zero, then for $R$ large enough there exist as much
non-trivial profiles as wanted. The critical value of $\alpha$ is obtained by
computing the unstable direction of $0$ in the phase plane. This direction is
given by the eigenvector $(1,\sqrt{ a\lambda})$ related to the
eigenvalue $\sqrt{ a\lambda}$ of the matrix 
$\left(\begin{array}{cc} 0&1\\ a\lambda & 0 \end{array}\right)$,
which corresponds to the linearization of the flow at the equilibrium
point $(0,0)$. This shows that $\DD$ intersects the homoclinic orbit to zero if
and only if $\alpha< a\lambda$.

The fact that there exists a profile $V_M$ larger than the other ones is
classical, see for example \cite{RBVL}. It is a consequence of Lemmas
\ref{lemme-compare} and \ref{lemme-gradient}: if $u(t)$ is a solution of
\eqref{eq-evo-prof-2} with $u(0)$ larger than any profile, then the
$\omega-$limit set $\omega(u(0))$ consists in a profile $V_M$ larger than any
profile. The fact that $V_M$ is even follows from the symmetry of the problem
since $V_M(-.)$ is also a profile. Notice that this symmetry corresponds to
the symmetry of the phase plane and to the fact that the trajectory $y\mapsto
\Phi_y(V_M(-R),V'_M(-R))$ intersect $\DD'$ for the first time at $y=R$.

%%%%%%%%%%%%%%%%%%%%%%%%%%%%%%%%%%%%%%%%%%%%%%%%%%%%%%%%%%%%%%%%

\subsection{A graphic criterium for stability of the
profiles}\label{subsection-exist-vp}

Let $V$ be a profile solution of \eqref{eq-profil}. We consider the
linearization $L_V$ of the flow along the profile $V$. We recall that
$L_V$ is given by \eqref{eq-LV}. The operator $L_V$ is self-adjoint and is a
relatively compact disturbance of $\Delta - \alpha$ and therefore has the
same essential spectrum, see \cite{Henry}. This yields the following
characterisation of the stability of a profile $V$, see \cite{Smoller}.
\begin{prop}
The spectrum of $L_V$ consists in $]-\infty ,-\alpha]$ and a finite number of
isolated real eigenvalues $\mu_i>-\alpha$ of finite multiplicity. As a
consequence, the profile $V$ is asymptotically stable if and only if $L_V$ has
no nonnegative eigenvalue.
\end{prop}

We are looking to nonnegative eigenvalues of $L_V$ that is to $\lambda\geq 0$
such that there exists $\varphi\in H^2(\Rm)$ such that $\varphi\not\equiv 0$ and
\begin{equation}\label{eq-vp-2}
\varphi''(y)+\left(f'(V(y))
\1_{|y|\leq R}-\alpha\1_{|y|>R}\right)\varphi(y)=\mu \varphi(y)~.
\end{equation}
By the same argument used for $V$ in Section
\ref{sect-para}, $\varphi$ is explicit outside $(-R,R)$
and it is equivalent to look to $\varphi$ satisfying 
\begin{equation}\label{eq-vp}
\left\{\begin{array}{ll} \varphi''(y)+f'(V(y))\varphi(y)=\mu \varphi(y) &
~~y\in(-R,R) \\ \varphi'(-R)=\sqrt{\mu+\alpha} \varphi(-R) & \\ 
\varphi'(R)=-\sqrt{\mu+\alpha} \varphi(R) & \end{array}\right.
\end{equation}
In this section, we adapt the standard Sturm-Liouville theory to the eigenvalue
problem \eqref{eq-vp}. Keeping in mind the geometric interpretation of
Sturm-Liouville arguments, we will obtain a graphic criterium to count
the number of nonnegative eigenvalues of $L_V$. The idea of using such graphic
arguments for studying the one-dimensional parabolic equation is not new, see
 \cite{Fusco-Rocha}. In fact, one
can understand the whole dynamics of a one-dimensional parabolic
equation by similar arguments as shown in \cite{Fiedler-Rocha}. See for example
\cite{Fiedler-Scheel} for a review on this subject.

For any $\mu\geq 0$ we define $\theta_\mu$ by
\begin{equation}\label{eq-theta}
\left\{\begin{array}{l} \theta_\mu'=-\sin^2
\theta_\mu+\left(\mu-f'(V(y))\right)\cos^2
\theta_\mu~~,~~~y\in(-R,R)\\
\theta_\mu(-R)=\arctan (\sqrt{\mu+\alpha}) \end{array}\right.
\end{equation}
and we set for $\mu>-\alpha$
\begin{equation}\label{eq-h}
h(\mu)=\theta_\mu(R)+\arctan(\sqrt{\mu+\alpha})~ .
\end{equation}
\begin{lemma}\label{lemme-1}
A nonnegative number $\mu$ is a nonnegative eigenvalue of $L_V$ if and only if
$h(\mu)$ belongs to $\pi\Zm$.
\end{lemma}
\begin{demo}
Let $\mu\geq 0$ and $\varphi\in\CC^2([-R,R])$, $\varphi\not\equiv
0$, satisfying \eqref{eq-vp}. We introduce two  functions of class
$\CC^2$, $\rho>0$ and $\theta\in\Rm$, such that
$(\varphi,\varphi')(y)=\rho(y)(\cos \theta(y), \sin
\theta(y))$. Notice that these functions exist and are regular since
the vector $(\varphi,\varphi')$ cannot vanish because $\varphi$ is a
non-trivial solution of a linear second-order differential
equation. Up to the change of the sign of $\varphi$ and to the
addition of a constant in $\pi\Zm$ to $\theta$,
$\varphi'(-R)=\sqrt{\mu+\alpha} \varphi(-R)$ implies
$\theta(-R)=\arctan (\sqrt{\mu+\alpha})$. Moreover, it follows
from \eqref{eq-vp} that, for $y\in (-R,R)$, 
\begin{align}
\theta'&=-\sin^2 \theta +\left(\mu-f'(V(y))\right)\cos^2 \theta \nonumber\\
\rho'&=\rho \sin \theta\cos\theta\left(1+\mu-f'(V(y))\right) \label{eq-lemme-1}
\end{align}
Hence, $\theta=\theta_\mu$ and $\varphi'(R)=-\sqrt{\mu+\alpha}
\varphi(R)$ is equivalent to
$\theta_\mu=-\arctan(\sqrt{\mu+\alpha})+k\pi$ with
$k\in\Zm$.\\ 
To show the other implication, it is sufficient to follow the previous
arguments in the opposite way: we start with a function
$\theta_\mu$, we construct $\rho$ by \eqref{eq-lemme-1} and we
check that $\varphi=\rho\cos\theta_\mu$ is a solution of
\eqref{eq-vp}. \end{demo} 

\begin{lemma}\label{lemme-2}
Let $V$ be a profile solution of \eqref{eq-profil}.\\ 
i) If $\mu>\|f'(V)\|_\infty$ then, for all $y\in [-R,R]$,
$\theta_\mu(y)\in(0,\frac \pi 2)$ and $h(\mu)\in (0,\pi)$.\\ 
ii) If $\mu_1>\mu_2$ then, for all $y\in [-R,R]$,
$\theta_{\mu_1}(y)>\theta_{\mu_2}(y)$ and $h$ is an increasing
function. 
\end{lemma}
\begin{demo}
We prove assertion i) by contradiction. Let $y_0=\inf \{y\in [-R,R]/
\theta_\mu(y)\not\in (0,\pi/2)\}$. Notice that $y_0>-R$ since
$\theta_\mu(-R)\in (0,\pi /2)$. If $\theta_\mu(y_0)=0$, then
$\theta'_\mu(y_0)=\mu-f'(V(y_0))>0$ which is absurd since
$\theta_\mu(y)>0$ for $y<y_0$ by definition of $y_0$. The
contradiction is similar if $\theta_\mu(y_0)=\pi/2$.\\ 
The proof of ii) is very similar: let $y_0=\inf \{y\in [-R,R]/
\theta_{\mu_1}(y)\leq \theta_{\mu_2}(y)\}$. We have
$\theta_{\mu_1}(y_0)=\theta_{\mu_2}(y_0)$. If
$\theta_{\mu_1}(y_0)\not\in \pi/2 + \pi\Zm$, then
$\theta'_{\mu_1}(y_0)-\theta'_{\mu_2}(y_0)=(\mu_1-\mu_2)\cos^2
\theta_{\mu_1}(y_0)>0$ and it contradicts the definition of
$y_0$. If $\theta_{\mu_i}(y_0)\in \pi/2 + \pi\Zm$, then
$\theta'_{\mu_1}(y_0)=\theta'_{\mu_2}(y_0)$, 
$\theta''_{\mu_1}(y_0)=\theta''_{\mu_2}(y_0)$ and
$\theta'''_{\mu_1}(y_0)=2+\mu_1-f'(V(y_0))>\theta'''_{\mu_2}(y_0)$. In
both cases, we obtain the desired contradiction. \end{demo} 

As a consequence of both preceding lemmas, we obtain the following
Sturm-Liouville-type result. 
\begin{prop}\label{prop-Sturm-Liouville}
Let $V$ be a profile solution of \eqref{eq-profil}. Let $L_V$ be the
operator defined by \eqref{eq-LV} and let $h$ be the function defined
by \eqref{eq-h}.\\ 
i) $k\in\Nm$ is such that $h(0)\in (-k\pi,(1-k)\pi]$ if and only if
  the operator $L_V$ has $k$ nonnegative eigenvalues.\\ 
ii) if $L_V$ has a nonnegative eigenvalue, there exists a positive
eigenfunction associated to the largest eigenvalue of $L_V$.\\
iii) the nonnegative eigenvalues of $L_V$ are simple.\\ 
\end{prop}
\begin{demo}
By Lemma \ref{lemme-1}, the
nonnegative eigenvalues correspond to the values $\mu$ for which
$h(\mu)\in \Zm$. Since, as shown in Lemma \ref{lemme-2}, $\mu\mapsto h(\mu)$ is
increasing and belongs to $(0,\pi)$ for large $\mu$, then Assertion i) follows
from the continuity of $h$. Moreover, if it exists,
the largest eigenvalue $\mu$ is necessary such that $h(\mu)=0$.
Thus $\theta_\mu(R)\in (-\pi/2,\pi/2)$. Since
$\theta_\mu(-R)\in(-\pi/2,\pi/2)$ and since $\theta_\mu'(y)<0$
as soon as $\theta_\mu(y)=\pm\pi/2$, we must have
$\theta_\mu(y)\in (-\pi/2,\pi/2)$ for all $y\in
[-R,R]$. Following the arguments in the proof of Lemma
\ref{lemme-1}, we construct an associated eigenfunction $\varphi=\rho
\cos \theta_\mu$ which is positive on $[-R,R]$. It is positive on
$\Rm$ by the extension imposed by \eqref{eq-vp-2}. Finally, the
eigenspace of a nonnegative eigenvalue is one-dimensional since any
eigenfunction must satisfy \eqref{eq-vp}. Since $L_V$ is
self-adjoint, any nonnegative eigenvalue is simple.\end{demo}

The interest of the first assertion of Proposition
\ref{prop-Sturm-Liouville} is to have a graphic
interpretation. Indeed, let $\varphi$ be the solution of
$(\varphi,\varphi')(-R)=(1,\sqrt{\alpha})$ and
$\varphi''+f'(V)\varphi=0$. In other words, $y\mapsto (\varphi,\varphi')(y)$ is
a trajectory of the linearization of the flow $D_V\Phi_y$ along the
profile $V$ such that $(\varphi,\varphi')(-R)$ is in the tangent space of $\DD$
at $(V,V')(-R)$. For each $y\in (-R,R)$, $(\varphi,\varphi')(y)$ is therefore
in the tangent space of $\Phi_y\DD$ at $(V,V')(y)$. By the same arguments as in
the proof of Lemma \ref{lemme-1}, $\theta_0(y)$ defined by \eqref{eq-theta} is
exactly the argument of the vector $(\varphi,\varphi')(y)$. 
Therefore, one can observe the angle $\theta_0(y)$ by
looking at the argument of the vector tangent to $\Phi_y\DD$ at
$(V,V')(y)$. Thus, we can graphically interpret Assertion i) of Proposition
\ref{prop-Sturm-Liouville} in the following way. 
\begin{graph-crit}\label{gc1}
Let $V$ be a profile, solution of \eqref{eq-profil} and let $L_V$ be
the operator defined by \eqref{eq-LV}. Let $y\in[-R,R] \longmapsto
A_V(y)\in\Rm^2$ be a continuous function such that for all $y\in
[-R,R]$, $A_V(y)$ is a unit
vector tangent to $\Phi_y\DD$ at $(V,V')(y)$. Then, the number of
times $A_V(y)$ crosses the direction of $\DD'$ when $y$ describes
$[-R,R]$ is exactly the number of nonnegative eigenvalues of
$L_V$. The crossings have to be counted in a algebraic way: positively
in the clockwise sense and negatively otherwise.
\end{graph-crit}

\begin{figure}[htp]
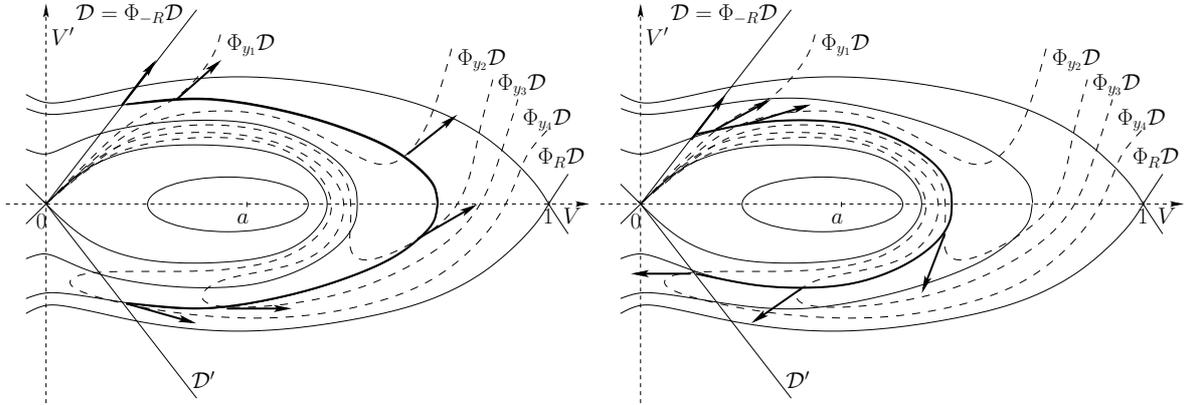

\begin{center}
\resizebox{0.48\textwidth}{!}{\input{count-vp-1.pstex_t}}
\resizebox{0.48\textwidth}{!}{\input{count-vp-2.pstex_t}}
%{\epsfig{file=count-ev-2.eps, width=7.8cm}}
\end{center}
\caption{\it an example of application of the graphic criterium \ref{gc1}. One
counts along a profile $V$ the number of crossings between the tangent vector to
the curve $\Phi_y\DD$ and the direction of $\DD'$. This gives the number of
nonnegative eigenvalues of $L_V$. Left, the largest profile has no nonnegative
eigenvalue. Right, the other non-trivial profile has one nonnegative
  eigenvalue.} 
\label{fig-vp1}
\end{figure}

The previous graphic criterium is illustrated in Figure \ref{fig-vp1}.
This criterium requires the knowledge of the curve
$\Phi_y\DD$ for each $y\in (-R,R)$. Therefore, it is difficult to apply as it
stands. Fortunately, there is a way to count the number of
times $A_V(y)$ crosses the direction of $\DD'$ with the knowledge of the curve
$\Phi_R\DD$ only. Indeed, in the particular case of this article, the phase
plane associated to $f$ is such that the curve $\Phi_R\DD$ is above $\DD'$ in
$\{(u,v)\in\Rm^2,~u\geq 1\}$ and below $\DD'$ in
$\{(u,v)\in\Rm^2,~u<0\}$. For this reason, using homotopy and topological
arguments, the vectors $A_0(y)$ and $A_{V_M}(y)$ corresponding to the extremal
profiles $0$ and $V_M$ must have a trivial number of crossings with the
direction of $\DD'$. Moreover, for another profile $V$, we can compute the
number of crossings of $A_V(y)$ with $\DD'$ as follows.
\begin{graph-crit}\label{gc2}
Let $V$ be a given profile. We follow the curve $\Phi_R\DD$ from $0$ (or from
$V_M$) to $V$. We count in a algebraic way the number of times the unit
tangent vector to $\Phi_R\DD$ crosses the direction of $\DD'$ during this
course. The resulting number is exactly the number of nonnegative eigenvalues
corresponding to the profile $V$.
\end{graph-crit}

The equivalence of the graphic criteria \ref{gc1} and \ref{gc2} can be easily
``seen on the figure''. However, the rigourous proof uses homotopies and the
topology of the plane, similarly to Jordan
theorem. This proof is not the subject of this paper and will be omitted.

\begin{figure}[htp]
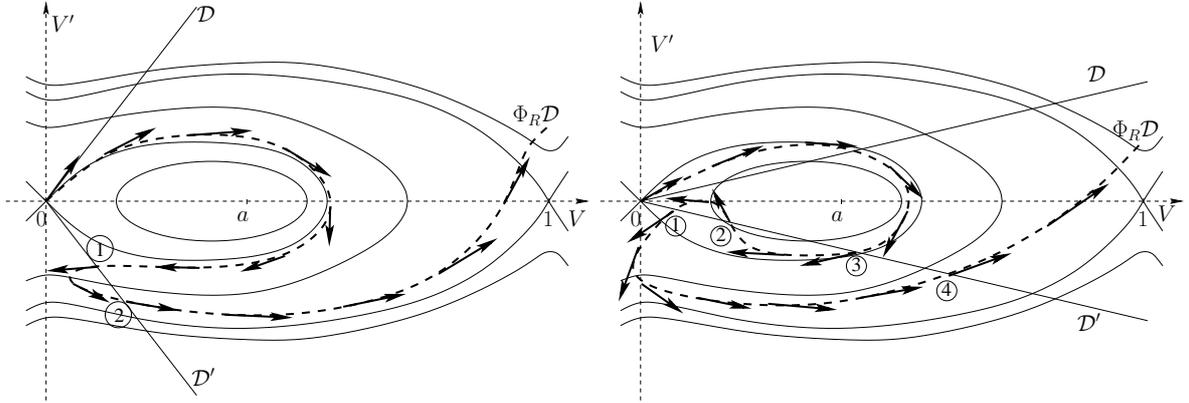

\begin{center}
\resizebox{0.48\textwidth}{!}{\input{count-vp-3.pstex_t}}
\resizebox{0.48\textwidth}{!}{\input{count-vp-4.pstex_t}}
%{\epsfig{file=count-ev-1.eps, width=7.8cm}} \hfill {\epsfig{file=count-ev-3.eps,width=7.8cm}}
\end{center}
\caption{\it examples of application of the graphic criterium \ref{gc2}. One
follows the curve $\Phi_R\DD$ from $0$ to a profile $V$ and counts
in a algebraic way the number of crossings between the tangent vector to
the curve and the direction of $\DD'$. One then obtains the number of
nonnegative eigenvalues of $L_V$. At the left, the largest profile (2)
has no nonnegative eigenvalue whereas the other non-trivial profile
(1) has one nonnegative eigenvalue. At the right, the four non-trivial
profiles have respectively one, two, one and zero nonnegative
eigenvalues.}
\label{fig-vp2}
\end{figure}

We can apply the above graphic criterium to the particular phase plane of this
article. Because the flow of $V''+f(V)$ is turning around the
equilibrium $(a,0)$, every profile $V$, which is not one of the extremal
profiles $0$ or $V_M$, is unstable, see Figure \ref{fig-vp2}.

%%%%%%%%%%%%%%%%%%%%%%%%%%%%%%%%%%%%%%%%%%%%%%%%%%%%%%%%%%%%%%%%
%%%%%%%%%%%%%%%%%%%%%%%%%%%%%%%%%%%%%%%%%%%%%%%%%%%%%%%%%%%%%%%%
%%%%%%%%%%%%%%%%%%%%%%%%%%%%%%%%%%%%%%%%%%%%%%%%%%%%%%%%%%%%%%%%

\section{Energy of the profiles}\label{sect-energy}
We recall that the relevant energy for the profiles is the function $E$
defined by \eqref{eq-energy}. We also recall that the energy $E(V)$ of a
profile $V$ is equal to the energy $\Em(V_{|(-R,R)})$ of the restriction of
$V$ to $(-R,R)$, see Section \ref{sect-para}. Let $R_0$ be the critical
thickness introduced in Proposition \ref{prop-profiles}. We know that for
$R>R_0$ there exists a unique
non-trivial profile $V_M$ which is stable and that $V_M$ is the
largest profile. This section is devoted to the following result.
\begin{prop}\label{prop-energie}
$~$\\
i) If $V$ is an unstable profile, then $E(V)$ is larger than $E(0)=0$
and $E(V_M)$.\\ 
ii) There exists $R_2\geq R_0$ such that on $(R_0,R_2)$ (this
set being possibly empty, for example if $\alpha\geq
 a\lambda$) the energy $E(V_M)$ of the largest profile is an
increasing function of the radius $R$ and on $(R_2,+\infty)$, $E(V_M)$
is a decreasing function of the radius $R$.\\ 
iii) $E(V_M)$ is positive for $R>R_0$ close to $R_0$ and converges to
$-\infty$ when $R$ goes to $+\infty$.\\
As a consequence, there exists a
radius $R_1>R_0$ such that $E(V_M)$ is positive for $R_0 < R \leq R_1$ and
negative for $R>R_1$. 
\end{prop}

The proof of Proposition \ref{prop-energie} relies on technics which are more
general than the framework of this paper. However, to obtain Proposition
\ref{prop-energie} we have to use properties of the phase plane corresponding to
\eqref{eq-profile-2}. These particular properties are observed on the phase
plane and on the numerical simulations. 

%%%%%%%%%%%%%%%%%%%%%%%%%%%%%%%%%%%%%%%%%%%%%%%%%%%%%%%%%%%%%%%%

\subsection{Energy of the unstable profiles}\label{subsect-instable}
We prove here the first assertion of Proposition
\ref{prop-energie}. The proof relies on the following lemma. 
\begin{lemma}\label{lemme-instable}
Let $V$ be an unstable profile. Let $\mu>0$ be the first
eigenvalue of $L_V$ and let $\varphi$ be a positive associated
eigenfunction, the existence of which is stated in Proposition
\ref{prop-Sturm-Liouville}. Then, there exist two globally bounded
solutions $u_-(t)$ and $u_+(t)$ of \eqref{eq-evo-prof-2} with the following
prescribed asymptotic behaviour in $H^1(-R,R)$: 
\begin{equation}\label{eq-asympt}
u_\pm(t)=V \pm e^{\mu t}\varphi + o\left(e^{\mu
  t}\right)~~~\text{when }t\longrightarrow -\infty~.
\end{equation}
As a consequence, for all $(y,t)\in(-R,R)\times \Rm$, $u_-(y,t)\leq V(y) \leq
u_+(y,t)$.
\end{lemma}
\begin{demo}
The lemma is a classical application of the theory of stable and
unstable manifolds near an equilibrium point. We refer to \cite{CHT}
and \cite{Henry}. Indeed, we can split the spectrum
of the linearisation $L_V$ into $\{\mu\}$ and the spectrum contained in
the half-plane $\{z\in\Cm,~ \Re(z) \leq \mu-\varepsilon\}$ with
$\varepsilon>0$ small enough. As shown in Proposition
\ref{prop-Sturm-Liouville}, $\mu$ is simple and admits a positive eigenfunction
$\varphi$. We defined the strongly unstable set by
\begin{align*}
W^{uu}(V)=\{ & u_0\in H^1(-R,R),~\exists~ u(t)\text{ global solution of
\eqref{eq-evo-prof-2} such that}\\
& \lim_{t\rightarrow -\infty} e^{(-\mu+\varepsilon/2)t}\|u(t)-V\|_{H^1(-R,R)} =
0~ \}~.
\end{align*}
The theory of invariant manifolds near an equilibrium point shows that
$W^{uu}(V)$ is an invariant one-dimensional manifold which is tangent
at $V$ to the line $V+\Rm \varphi$. The manifold $W^{uu}(V)$ consists in $V$
and two globally defined trajectories $u_+$ and $u_-$ satisfying the asymptotic
behaviour \eqref{eq-asympt}. Of course, the last assertion is a direct
consequence of this asymptotic behaviour, of the positivity of $\varphi$ and of
Lemma \ref{lemme-compare}. \end{demo} 

The first assertion of Proposition \ref{prop-energie} is deduced from Lemma
\ref{lemme-instable} as follows.
Let $V$ be an unstable profile. We know from Section \ref{sect-exist}
that $V$ lies between $0$ and the largest profile $V_M$, both being
stable. Let $u_-(t)$ be the solution given by Lemma
\ref{lemme-instable} for the profile $V$. By Lemma
\ref{lemme-compacite}, the $\omega-$limit set of $u_-$ is non-empty
and consists in profiles. We wonder which profiles may belong to this
$\omega-$limit set. Assume that there is an unstable profile $\tilde
V$ in $\omega(u_-)$. This profile must satisfy $\tilde V\leq V$ since
$u_-(t)\leq V$ for all $t$. Moreover, by Lemma \ref{lemme-gradient},
$\tilde V\not \equiv V$ and thus, by uniqueness of ODE solutions,
$\tilde V(y)<V(y)$ for all $y\in\Rm$. Let $\tilde u_+(t)$ be the
solution given by Lemma \ref{lemme-instable} for the profile $\tilde
V$. For $t$ close to $-\infty$, $\tilde
u_+(t)<u_-(t)$. Using Lemma
\ref{lemme-compare}, $\tilde V \leq \tilde u_+(t) \leq  u_-(t) \leq V$
for all times. Since Lemma \ref{lemme-gradient} prevents $\tilde u_+(t)$ to go
back to $\tilde V$, $u_-(t)$ cannot converge to $\tilde V$, which yields a
contradiction. Therefore, since $\omega(u_-)$
is non-empty, the only possibility is $\omega(u_-)=\{0\}$. Lemma
\ref{lemme-gradient} then shows that $\Em(V)>\Em(0)=0$. 

By using $u_+(t)$, we prove similarly that $\Em(V)>\Em(V_M)$. Notice that
the above arguments also prove that every unstable profile intersects
the other unstable ones as illustrated in Figure \ref{fig-R3}. 

%%%%%%%%%%%%%%%%%%%%%%%%%%%%%%%%%%%%%%%%%%%%%%%%%%%%%%%%%%%%%%%%

\subsection{Energy of the stable profile}\label{subsect-stable}
The proof of the second and third assertions of Proposition
\ref{prop-energie} is splitted in several parts. We use in this
section the notations of Section \ref{subsection-exist-profile} with
obvious changes  to add the dependence of the different objects with
respect to $R$.\\

\begin{figure}[htp]
\begin{center}
\resizebox{0.6\textwidth}{!}{\input{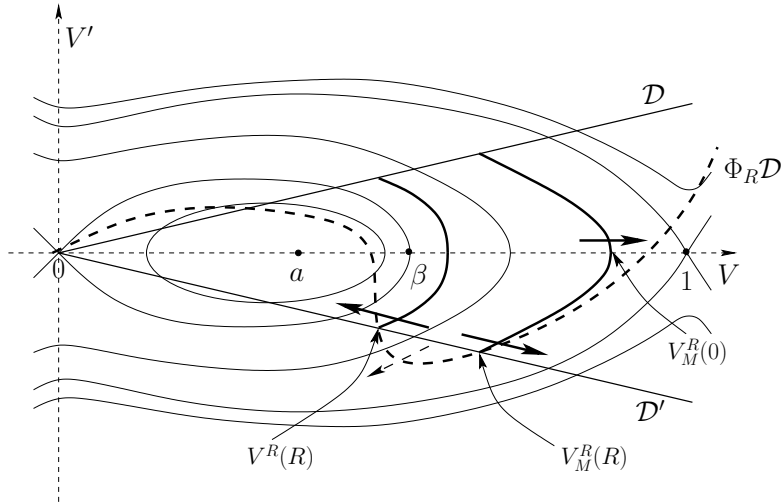}}
%\resizebox{0.7\textwidth}{!}{\input{fig-energy.pstex_t}}
\end{center}
\caption{\it when $R=R_0$, a saddle-node bifurcation creates a stable profile
$V^R_M$ (the largest profile) and an unstable profile $V^R$. When $R>R_0$
increases, the largest profile $V^R_M$ is increasing. When $R\rightarrow
+\infty$ the largest profile $V^R_M$ piles up to the singularity $(1,0)$ of the
phase plane.}
\label{fig-energy}
\end{figure}

\noindent $\bullet$ {\bf $R\mapsto V^R_M(0)$ is increasing.}\\
To study the function $R\mapsto V^R_M$, we invoke the graphic
interpretation of the profiles developped in Section
\ref{subsection-exist-profile}. By symmetry of the largest profile, $V^R_M(0)$
is the intersection of the trajectory $y\mapsto \Phi_y (V^R_M(-R),\alpha
V^R_M(-R))$ with the horizontal axis. Therefore, $R\mapsto V^R_M(0)$ is
increasing if and only if $R\mapsto V^R_M(-R)$ is increasing. For $R$ given,
$\Phi_R(V^R_M(-R),\alpha
V^R_M(-R))$ belongs to $\DD'$ and at this point of the phase plane,
the flow $\Phi_y$ crosses $\DD'$ transversally from the top to the
bottom. Thus, for small $\varepsilon>0$,
$\Phi_{R+\varepsilon}(V^R_M(-R),\alpha V^R_M(-R))$ is below $\DD'$ and
by continuity, there exists an intersection of
$\Phi_{R+\varepsilon}\DD$ and $\DD'$ at the right of
$\Phi_{R+\varepsilon}(V^R_M(-R),\alpha V^R_M(-R))$. This means that
$V^{R+\varepsilon}_M(-R-\varepsilon)>V^R_M(-R)$ and thus $R\mapsto
V^R_M(0)$ is increasing. See Figure \ref{fig-energy}.

By a standard implicit functions argument, one shows that $R\mapsto V^R_M$ is
continuous and differentiable as soon as $\Phi_{R}\DD$ intersects
$\DD'$ transversally at $(V^R_M(R),\alpha V^R_M(R))$. We admit that it is
always the case when $R>R_0$.

\noindent$\bullet$ {\bf $E(V^R_M)$ converges to $-\infty$ when $R$ goes
  to $+\infty$.}\\ 
This fact has already been proved in \cite{GC} by the first author. We only
give here brief arguments. When $R$ increases, the trajectory $y\in(-R,R)\mapsto
\Phi_y(V^R_M(-R),\alpha V^R_M(-R))$ goes to the right and $V^R_M(0)$
converges to $1$. For large $R$, the trajectory $y\in(-R,R)\mapsto
\Phi_y (V^R_M(-R),\alpha V^R_M(-R))$ spends a time $O(R)$ in a neighborhood of
$(1,0)$, where $(u,v)\mapsto \frac 12 |v|^2 - F(u)$ is negative, and a bounded
time outside. Therefore $E(V^R_M)=\Em(V^R_M)$ goes to $-\infty$ when $R$ goes to
$+\infty$.\\

\noindent$\bullet$ {\bf There exists $\beta\in(a,1)$ such that $R\mapsto
E_R(V^R_M)$ is increasing if $V^R_M(0)<\beta$ and decreasing if
$V^R_M(0)>\beta$}\\ 
We admit that for all $R>R_0$, $R\mapsto V^R_M(y)$ is differentiable. Let
$E_\rho$ be the energy functional defined by \eqref{eq-energy} with
$R=\rho$. Let $V$ be a given fixed profile. We consider the variation of $\rho
\mapsto E_\rho(V)$ at $\rho=R$.
$$\partial_\rho E_\rho(V)_{|\rho=R} = \left(-F(V(R))-\frac
\alpha 2 \left|V(R)\right|^2 \right)+\left(-F(V(-R))-\frac
\alpha 2 \left|V(-R)\right|^2\right)
$$
Since $V(\pm R)=\pm\sqrt\alpha\;V'(\pm R)$ and since $V''(y)+f(V(y))=0$,
$y\mapsto \frac 12\left|\partial_y V(y)\right|^2 +F(V(y))$ is constant and   
$$\partial_\rho E_\rho(V)=-2 F(V(0))~.$$
Since $V$ is an equilibrium of \eqref{eq-evo-prof}, it is a critical
point for the energy $E_R$. In other words, for any $W\in H^1(\Rm)$,   
$$\partial_h E_R(V+ h W)=0~.$$
Henceforth, if $R\mapsto V^R_M$ is differentiable, then combining both
derivatives yields 
$$\partial_R \left(E_R(V^R_M)\right)= -2 F(V^R_M(0))~.$$
It remains to notice that $R\mapsto V^R_M(0)$ is increasing and
converges to $1$ when $R$ goes to $+\infty$. Moreover, there exists
$\beta>0$ such that $F$ is negative on $(0,\beta)$ and positive on
$(\beta,1]$. In fact, $(\beta,0)$ is the point of the phase plane where the
homoclinic orbit to 0 intersects the horizontal axis, see Figure
\ref{fig-energy}.

\noindent$\bullet$ {\bf$E_R(V^R_M)$ is positive for $R$ close to $R_0$.}\\
When $R$ passes the value $R_0$, the creation of both non trivial profiles
occurs through a saddle-node bifurcation: two equilibrium states appear at the
same point, one stable and one unstable. The stable profile is $V^R_M$, the
largest one, and we denote by $V^R$ the unstable profile which lies between $0$
and $V^R_M$. From now on, we work in $(-R,R)$ by using the analogy
presented in Section \ref{sect-para}. We want to show that $\Em(V^R_M)$ is
positive for $R$ close to $R_0$. When $R$ decreases to $R_0$: the profiles
$V_M^R$ and $V^R$ collide. As $\Em(V^R)>0$ for $R>R_0$ as shown in Section
\ref{subsect-instable}, we must have at the limit
$\Em(V^{R_0})=\Em(V^{R_0}_M)\geq 0$. Assume that $\Em(V^{R_0}_M)=0$. Let
$R_n>R_0$ be a sequence of thickness decreasing to $R_0$ and let $V^{R_n}$ be
the associated unstable profiles. Let $u^n_-(t)$ be the sequence of solutions of
\eqref{eq-evo-prof-2} given by Lemma \ref{lemme-instable}
applied to $V^{R_n}$. We know by Section \ref{subsect-instable} that
$u^n_-(t)$ converges to $0$ when $t$ goes to $+\infty$. Let
$K=\lim_{n\rightarrow +\infty} \|V^{R_n}\|_{H^1(-R,R)}$. Notice that $K$ is
positive since for each $n$ there exists $y\in(-R,R)$ such that $V^{R_n}(y)>a$. 
We set $t_n$ to be a time such that $\|u^n_-(t_n)\|_{H^1(-R,R)}=K/2$. By a
compacity argument similar to Lemma \ref{lemme-compacite} ($R$ is moving, but
nothing singular happens), one can extract a subsequence such that
$u^{\varphi(n)}_-(t_{\varphi(n)})$ converges in $H^1(-R,R)$ to a
function $u^\infty$. Notice that by construction, $u^\infty$ is
neither $0$ nor $V^{R_0}$. The gradient structure of \eqref{eq-evo-prof-2}
shows that, for all $n$ and $t$, $0<\Em(u^n_-(t))<\Em(V^{R_n})$ and thus
$\Em(u^\infty)=0$. Let $u^\infty(t)$ be the solution of \eqref{eq-evo-prof-2}
for $R=R_0$, with initial data $u^\infty$. This solution is not a profile,
and so its energy $\Em(u(t))$ decreases and is negative for $t>0$. However,
$u^\infty(t)$
must converge to a profile when $t$ goes to $+\infty$ due to Lemma
\ref{lemme-gradient}. Since all the profiles at $R=R_0$ have an energy equal to
$0$ by assumption, this is impossible and we get a
contradiction. Therefore, $\Em(V^{R_0})=\Em(V^{R_0}_M)=E(V^{R_0}_M)$ must be
positive. 

\section{Discussion}\label{sect-discussion}

In this paper, we have proved that there exist two critical
thicknesses $R_1>R_0>0$ such that if $0<R<R_0$, there is no
non-trivial profile solution of \eqref{eq-profil}. If
$R_0<R$, there exist non-trivial profiles. One of them is larger than
every other one and is stable, whereas every other non-trivial profiles are
unstable. Finally, the energy of the unstable profiles is always
larger than the energy of the stable profiles. If $R_0<R<R_1$, the
energy of the largest profile is
larger than the energy of $0$, whereas it is smaller if $R>R_1$.
These results give us
informations on the propagation of travelling fronts solution of equation
\eqref{eq} as stated in Consequence \ref{th-2}.

We recall that in equation \eqref{eq}, $u$ represents the
depolarization of the brain so if $u(X)=0$ the brain is normally
polarized at the point $X$, and if $u(X)=1$ the brain is totally depolarized. A
depolarization wave in the brain corresponds to a travelling front
solution of \eqref{eq} where a non-trivial stable profile invades
the zero state. Thus, the above mathematical analysis of the stability and
energy of the asymptotic profiles of equation \eqref{eq-profil} yields
informations on the propagation of depolarization waves in the human
brain during stroke. If the width of the grey
matter 
is smaller than $R_1$ no depolarization wave can propagate through the
brain whereas if the width of the grey matter is larger than $R_1$, there exist
depolarization waves. This may explain why attempts to observe
depolarization waves in the human brain have received opposite
conclusions in different studies \cite{aitken, mayevsky, maclachlan, 
somjen, sramka}. In the human brain, the grey matter is particularly
thin and its thickness may vary a lot. The
difficulties to observe depolarization waves in the human brain could
thus be explained simply by the morphology of the brain. If the
initial depolarization takes place in a part of the brain where the
grey matter is large enough, then the depolarization can spread
through a part of the grey matter. But if the initial depolarization
takes place in a part of the brain where the grey matter is very thin
then no propagation will occur. Hence it is possible to observe
depolarization waves in the human brain \cite{aitken, mayevsky,
  maclachlan}, but they will not appear in all the cases \cite{somjen,
  sramka}. Moreover the depolarizations would not travel over large
distances because they will stop as soon as the grey matter becomes
too thin. For example, they will stop at the bottom of the large
sulkus, where the grey matter becomes thiner. This has already be
observed in the case of the migraine with aura. The aura may be due to
a depolarization wave \cite{Gorji, James} and for most of the patients the 
aura stops at the bottom of the Rolando sulkus.

This paper also gives informations on how the initial depolarization
is erased as illustrated in Figure \ref{fig-intro}. If the grey matter is very
thin, then the depolarization is
quickly absorbed uniformly in the excited area. If the grey
matter is a little bit larger, the depolarization occurs in a different way. The
depolarized area shrinks progressively while the cells stay
totally depolarized as long as they are in the middle of the excited area.
In this case, the repolarization is due to a travelling front where the
normally polarized state invades the depolarized state. The neurons in the
center of the excited area may stay depolarized a very long time, which may
cause local damages even if no spreading wave is observed. To our
knowledge this behaviour has never been observed in experiments.

In order to verify biologically these results and to validate the
model, it would be interesting 
to estimate numerically the values of the thresholds $R_0$ and $R_1$. They
depend on the values of the parameters of equation \eqref{eq} that are
hard to compute due to the few possibilities of quantitative mesures
in the brain. Moreover these parameter
values may vary a lot from one person to another and from a species to
another. Estimations of this parameter values for the rodent and for
the human may give us more informations on the difficulties to observe
depolarization waves in the human brain. Indeed
even if the mecanisms are the same in the rodent brain and in the
human brain, the densities of cells are totally different and this 
must influence the values of the parameters of this model and thus of
the thresholds. For the moment we are just trying to understand
numerically the influence of each parameters on these thresholds.

\vspace{5mm}

\noindent {\bf Acknowledgements: } the authors would like to thank Thierry
Gallay for several fruitful discussions.

%%%%%%%%%%%%%%%%%%%%%%%%%%%%%%%%%%%%%%%%%%%%%%%%%%%%%%%%%%%%%%%%%

%%%%%%%%%%%%%%%%%%%%%%%%%%%%%%%%%%%%%%%%%%%%%%%%%%%%%%%%%%%%%%%%%
%%%%%%%%%%%%%%%%%%%%%%%%%%%%%%%%%%%%%%%%%%%%%%%%%%%%%%%%%%%%%%%%%

\end{document}